\documentclass[11pt]{article}

\usepackage{latexsym}
\usepackage{amssymb}

\newtheorem{thm}{Theorem}
\newtheorem{prop}{Proposition}
\newtheorem{cor}{Corollary}

\begin{document}
{
\begin{center}
{\Large\bf
On the truncated operator trigonometric moment problem.}
\end{center}
\begin{center}
{\bf S.M. Zagorodnyuk}
\end{center}

\section{Introduction.}

The truncated operator trigonometric moment problem
consists of finding a non-decreasing $[\mathcal{H}]$-valued function
$F(t)$, $t\in [0,2\pi]$, $F(0)=0$, which is
strongly left-continuous in $(0,2\pi]$ and such that
\begin{equation}
\label{f1_1}
\int_0^{2\pi} e^{int} dF(t) = S_n,\qquad n=0,1,...,d,
\end{equation}
where $\{ S_n \}_{n=0}^d$ is a prescribed sequence of bounded operators on $\mathcal{H}$ (moments).
Here $\mathcal{H}$ is a fixed Hilbert space and $d\in \mathbb{Z}_+$ is a fixed number.
The operator Stieltjes integrals in~(\ref{f1_1}) are understood as limits of the corresponding
integral sums in the strong operator topology.

\noindent
Set $S_{-k} = S_k^*$, $k=1,2,...,d$.
The following condition:
\begin{equation}
\label{f1_2}
\sum_{k,l=0}^d (S_{l-k} h_l, h_k)_{\mathcal{H}} \geq 0,
\end{equation}
where $\{ h_k \}_0^d$ are arbitrary elements of $\mathcal{H}$,
is necessary and sufficient for the solvability of the moment problem~(\ref{f1_1}) (e.g.~\cite{cit_1000_A}).
The  solvable moment problem~(\ref{f1_1}) is said to be \textit{determinate} if it has a unique solution and \textit{indeterminate}
in the opposite case.
The truncated operator trigonometric moment problem was studied in papers~\cite{cit_1000_A}, \cite{cit_5000_IT} (in slightly different
statements).
The conditions of the solvability were obtained in~\cite{cit_1000_A}. In the case of the strict positivity of
the corresponding Toeplitz operator, all solutions to the moment problem were described in~\cite{cit_5000_IT}.
For the history of scalar and matrix truncated trigonometric moment problems we refer to~\cite{cit_10000_Z}, \cite{cit_15000_Z}.

Our aim here is to describe all solutions to the solvable moment problem~(\ref{f1_1}) without any additional conditions.
For this purpose we shall develop the operator approach of
Sz\"okefalvi-Nagy and Koranyi in~\cite{cit_7000_SK},~\cite{cit_9000_SK}.
Also our approach is close to the approach of Krein and Krasnoselskii in~\cite{cit_6000_KK}.
All solutions of the moment problem
are described by a Nevanlinna-type parameterization. In the case of moments acting in a separable Hilbert space, the matrices of
the operator coefficients in the Nevanlinna-type formula are calculated by the prescribed moments.
Conditions for the determinacy of the moment problem are given.

\noindent
\textbf{Notations.}
As usual, we denote by $\mathbb{R}, \mathbb{C}, \mathbb{N}, \mathbb{Z}, \mathbb{Z}_+$,
the sets of real numbers, complex numbers, positive integers, integers and non-negative integers,
respectively; $\mathbb{D} = \{ z\in \mathbb{C}:\ |z|<1 \}$; 
$\mathbb{T} = \{ z\in \mathbb{C}:\ |z|=1 \}$; $\mathbb{T}_e = \{ z\in \mathbb{C}:\ |z|\not= 1 \}$.
By $k\in\overline{0,\rho}$ we mean that $k\in \mathbb{Z}_+:\ 0\leq k\leq \rho$ if $\rho\in\mathbb{Z}_+$,
and $k\in \mathbb{Z}_+$, if $\rho = +\infty$.

\noindent
In this paper Hilbert spaces are not necessarily separable, operators in them are
supposed to be linear.

\noindent
If H is a Hilbert space then $(\cdot,\cdot)_H$ and $\| \cdot \|_H$ mean
the scalar product and the norm in $H$, respectively.
Indices may be omitted in obvious cases.
For a linear operator $A$ in $H$, we denote by $D(A)$
its  domain, by $R(A)$ its range, and $A^*$ means the adjoint operator
if it exists. If $A$ is invertible then $A^{-1}$ means its
inverse. $\overline{A}$ means the closure of the operator, if the
operator is closable. If $A$ is bounded then $\| A \|$ denotes its
norm.
For a set $M\subseteq H$
we denote by $\overline{M}$ the closure of $M$ in the norm of $H$.
For an arbitrary set of elements $\{ x_n \}_{n\in I}$ in
$H$, we denote by $\mathop{\rm Lin}\nolimits\{ x_n \}_{n\in I}$
the set of all linear combinations of elements $x_n$,
and $\mathop{\rm span}\nolimits\{ x_n \}_{n\in I}
:= \overline{ \mathop{\rm Lin}\nolimits\{ x_n \}_{n\in I} }$.
Here $I$ is an arbitrary set of indices.
By $E_H$ we denote the identity operator in $H$, i.e. $E_H x = x$,
$x\in H$. In obvious cases we may omit the index $H$. If $H_1$ is a subspace of $H$, then $P_{H_1} =
P_{H_1}^{H}$ is an operator of the orthogonal projection on $H_1$
in $H$.
By $[H]$ we denote a set of all bounded operators on $H$.
For a closed isometric operator $V$ in $H$ we denote:
$M_\zeta(V) = (E_H - \zeta V) D(V)$, $N_\zeta(V) = H\ominus M_\zeta(V)$, $\zeta\in \mathbb{C}$; $M_\infty(V)=R(V)$, $N_\infty(V)= H\ominus R(V)$.

\noindent
By $\mathcal{S}(D;N,N')$ we denote a class of all analytic in a domain $D\subseteq \mathbb{C}$
operator-valued functions $F(z)$, which values are linear non-expanding operators mapping the whole
$N$ into $N'$, where $N$ and $N'$ are some Hilbert spaces.

\section{The solvability and a description of solutions for the moment problem.}

Suppose that the moment problem~(\ref{f1_1}) is given and it is solvable. For arbitrary elements $\{ h_k \}_{0}^d$ of $\mathcal{H}$
we may write:
$$ \sum_{k,l=0}^d (S_{l-k} h_l, h_k)_{\mathcal{H}} =
\sum_{k,l=0}^d \left( \int_0^{2\pi} e^{i(l-k)t} dF(t) h_l, h_k \right)_{\mathcal{H}} = $$
$$ = \sum_{k,l=0}^d \left( \lim_{\delta\to +0} \sum_{r=0}^N  e^{i(l-k)t_r} (F(t_{r+1}) - F(t_{r})) h_l, h_k \right)_{\mathcal{H}} = $$
$$ = \lim_{\delta\to +0} \sum_{r=0}^N \left( 
(F(t_{r+1}) - F(t_{r})) \sum_{l=0}^d e^{i l t_r} h_l, \sum_{k=0}^d e^{i k t_r} h_k 
\right)_{\mathcal{H}} \geq 0. $$
Here $\delta$ is the diameter of a partition of $[0,2\pi]$ and $\{ t_r \}_0^N$ are points of this partition. Thus, condition~(\ref{f1_2}) is satisfied.

Conversely, suppose that the moment problem~(\ref{f1_1}) with $d\in N$ is given and condition~(\ref{f1_2}) is satisfied.
Like it was done in~\cite{cit_9000_SK} we consider abstract symbols $\varepsilon_j$, $j=0,1,...,d$, and form a formal sum $h$:
\begin{equation}
\label{f2_1}
h = \sum_{j=0}^d h_j \varepsilon_j, 
\end{equation}
where $h_j\in \mathcal{H}$.
If $\alpha\in \mathbb{C}$, then we set 
$\alpha h = \sum_{j=0}^d (\alpha h_j) \varepsilon_j$.
If
\begin{equation}
\label{f2_2}
g = \sum_{j=0}^d g_j \varepsilon_j, 
\end{equation}
where $g_j\in \mathcal{H}$, then
we set
$h+g =  \sum_{j=0}^d (h_j + g_j) \varepsilon_j$. 
A set of all formal sums of type~(\ref{f2_1}) becomes a complex linear vector space $\mathfrak{B}$.
Let $h,g\in \mathfrak{B}$ have the form as in~(\ref{f2_1}), (\ref{f2_2}).
Let
\begin{equation}
\label{f2_3}
\Phi(h,g) = \sum_{j,k=0}^d  ( S_{j-k} h_j, g_k )_{\mathcal{H}}. 
\end{equation}
The functional $\Phi$ is sesquilinear and it has the properties $\overline{\Phi(h,g)} = \Phi(g,h)$, $\Phi(h,h)\geq 0$.
If $\Phi(h-g,h-g)=0$, we put elements $h$ and $g$ to the same equivalence class denoted by $[h]$ or $[g]$.
A set of all equivalent classes we denote by $\mathfrak{L}$. By the completion of $\mathfrak{L}$ we obtain
a Hilbert space $H$.
Set
\begin{equation}
\label{f2_4}
x_{h,j}  := [ h \varepsilon_j ],\qquad h\in\mathcal{H},\ 0\leq j\leq d.
\end{equation}
Observe that
\begin{equation}
\label{f2_5}
( x_{h,j}, x_{g,k} )_H = ( S_{j-k} h, g )_{\mathcal{H}},\qquad h,g\in\mathcal{H},\ 0\leq j,k\leq d.
\end{equation}
Set 
$$ D_0 = \mathop{\rm Lin}\nolimits \{ x_{h,j} \}_{h\in \mathcal{H},\ 0\leq j\leq d-1} = $$
$$ = \{ x_{h_0,0} + x_{h_1,1} + ... + x_{h_{d-1},d-1} \}_{ h_0,h_1,...,h_{d-1}\in \mathcal{H} }. $$
Consider a linear operator $A_0$ with $D(A_0) = D_0$:
\begin{equation}
\label{f2_6}
A_0 \sum_{j=0}^{d-1} x_{h_j,j} = \sum_{j=0}^{d-1} x_{h_j,j+1},\qquad h_j\in \mathcal{H}.
\end{equation}
Let us check that $A_0$ is well-defined.
Suppose that an element $h\in D_0$ has two representations:
$$ h = \sum_{j=0}^{d-1} x_{h_j,j} = \sum_{j=0}^{d-1} x_{g_j,j},\qquad h_j,g_j\in \mathcal{H}. $$
Then
$$ \left\| 
\sum_{j=0}^{d-1} x_{h_j,j+1} - \sum_{j=0}^{d-1} x_{g_j,j+1}
\right\|_H
=
\left( \sum_{j=0}^{d-1} x_{h_j-g_j,j+1}, \sum_{k=0}^{d-1} x_{h_k-g_k,k+1} \right)_H = $$
$$ = \sum_{j,k=0}^{d-1} ( S_{j-k} (h_j-g_j), h_k-g_k )_{\mathcal{H}} =
\left( \sum_{j=0}^{d-1} x_{h_j-g_j,j}, \sum_{k=0}^{d-1} x_{h_k-g_k,k} \right)_H = 0. $$
Thus, $A_0$ is well-defined.
If
$$ h = \sum_{j=0}^{d-1} x_{h_j,j},\quad f = \sum_{k=0}^{d-1} x_{f_k,k},\qquad h_j,g_k\in \mathcal{H}, $$
then
$$ (A_0 h, A_0 f)_H = \sum_{j,k=0}^{d-1} (x_{h_j,j+1}, x_{f_k,k+1})_H = 
\sum_{j,k=0}^{d-1} ( S_{j-k} h_j, f_k )_{\mathcal{H}} = $$
$$ = \sum_{j,k=0}^{d-1} (x_{h_j,j}, x_{f_k,k})_H = (h,g)_H. $$
Therefore $A_0$ is isometric.
Set $A = \overline{A_0}$. 
By the induction argument it may be checked that
$$ x_{h,j} = A^j x_{h,0},\qquad 0\leq j\leq d. $$
Let $\widetilde A\supseteq A$ be a unitary operator in a Hilbert space
$\widetilde H\supseteq H$, and $\{ \widetilde E_t \}_{t\in [0,2\pi]}$ be its strongly left-continuous orthogonal resolution of the identity. 
We may write
$$ (S_j h,g)_{\mathcal{H}} = (x_{h,j}, x_{g,0})_H = (A^j x_{h,0},x_{g,0})_H = (\widetilde A^j x_{h,0},x_{g,0})_{\widetilde H} = $$
\begin{equation}
\label{f2_7}
= \left( \int_0^{2\pi} e^{itj} d \widetilde E_t x_{h,0},x_{g,0} \right)_{\widetilde H},\qquad
0\leq j\leq d,\ h,g\in \mathcal{H}. 
\end{equation}
Consider the following operator $I$: $\mathcal{H}\rightarrow H$:
\begin{equation}
\label{f2_8}
I h = x_{h,0},\qquad h\in \mathcal{H}.
\end{equation}
It is readily checked that $I$ is linear. Moreover, since
$$ \| Ih \|_H^2 = (x_{h,0}, x_{h,0})_H = (S_0 h,h)_{\mathcal{H}} \leq \| S_0 \| \| h \|_{\mathcal{H}}^2, $$
then $I$ is bounded.
By~(\ref{f2_7}) we may write
$$ (S_j h,g)_{\mathcal{H}} =
\left( P^{\widetilde H}_H \int_0^{2\pi} e^{itj} d \widetilde E_t Ih, Ig \right)_H =
\left( I^* P^{\widetilde H}_H \int_0^{2\pi} e^{itj} d \widetilde E_t Ih, g \right)_{\mathcal{H}} = $$
$$ = \left( \int_0^{2\pi} e^{itj} d \left( I^* P^{\widetilde H}_H \widetilde E_t I \right) h, g \right)_{\mathcal{H}},\qquad
0\leq j\leq d,\ h,g\in \mathcal{H}. $$
Therefore
\begin{equation}
\label{f2_9}
S_j = \int_0^{2\pi} e^{itj} d \left( I^* \mathbf{E}_t I \right),\qquad 0\leq j\leq d,
\end{equation}
where $\mathbf{E}_t$ is a strongly left-continuous spectral function of $A$ (corresponding to $\widetilde A$).
Thus,
\begin{equation}
\label{f2_10}
F(t) = I^* \mathbf{E}_t I,\qquad t\in [0,2\pi],
\end{equation}
is a solution of the moment problem~(\ref{f1_1}).
We conclude that each strongly left-continuous spectral function of $A$ generates a solution of
the moment problem~(\ref{f1_1}) by relation~(\ref{f2_10}).

Let $\widehat F(t)$ be an arbitrary solution of the moment problem~(\ref{f1_1}). We shall check that
$\widehat F(t)$ can be constructed by relation~(\ref{f2_10}).
By $C_{00} (\mathcal{H};[0,2\pi])$ we denote a set of all strongly continuous $\mathcal{H}$-valued functions $f(t)$, $t\in [0,2\pi]$,
which take their values in finite-dimentional subspaces of $\mathcal{H}$ (depending on $f$).
For arbitrary $f,g\in C_{00} (\mathcal{H};[0,2\pi])$ we set (see~\cite{cit_700_Ber})
\begin{equation}
\label{f2_11}
\Psi(f,g) = \int_0^{2\pi} \left(
d\widehat F(t) f(t), g(t)
\right)_{\mathcal{H}} := 
\lim_{\delta\to +0} \sum_{k=1}^N 
\left(
\widehat F(\Delta_k) f(t_k), g(t_k)
\right)_{\mathcal{H}},
\end{equation}
where $\delta$ is the diameter of a partition $\{ \Delta_k \}_{k=1}^N$, $\Delta_k = [a_k,b_k)$  of $[0,2\pi)$ and $t_k\in\Delta_k$.
Here, as usual, the limit does not depend on the choice of partitions and points $t_k$.
It is easy to see that in the case of $f,g\in C_{00} (\mathcal{H};[0,2\pi])$ the limit in~(\ref{f2_11}) exists and
reduces to a finite sum of scalar Stieltjes-type integrals.

Introducing classes of the equivalence with respect to $\Psi$ and by the completion we obtain 
a Hilbert space $L_2 = L_2(\mathcal{H}; [0,2\pi]; d\widehat F(t))$.
Consider two operator polynomials of the following form:
\begin{equation}
\label{f2_12}
p(t) = \sum_{j=0}^d e^{ijt} h_j,\quad q(t) = \sum_{k=0}^d e^{ikt} g_k,\qquad h_j,g_k\in \mathcal{H}.
\end{equation}
Since $p,q\in C_{00} (\mathcal{H};[0,2\pi])$ then the corresponding classes $[p],[q]$ belong to $L_2(\mathcal{H}; [0,2\pi]; d\widehat F(t))$.
As usual in such situations, we shall say that $p,q$ belong to $L_2(\mathcal{H}; [0,2\pi]; d\widehat F(t))$.
Then
$$ (p,q)_{L_2(\mathcal{H}; [0,2\pi]; d\widehat F(t))}
= \int_0^{2\pi} \left(
d\widehat F(t) p(t), q(t)
\right)_{\mathcal{H}}
= $$
$$ = \sum_{j,k=0}^d \int_0^{2\pi} e^{i(j-k)t} d \left( \widehat F(t) h_j, g_k \right)_{\mathcal{H}} = 
\sum_{j,k=0}^d ( S_{j-k} h_j, g_k )_{\mathcal{H}} = $$
\begin{equation}
\label{f2_13}
= \sum_{j,k=0}^d (x_{h_j,j}, x_{g_k,k})_H = \left( \sum_{j=0}^d x_{h_j,j}, \sum_{k=0}^d x_{g_k,k} \right)_H. 
\end{equation}
Denote by $P (\mathcal{H}; [0,2\pi]; d\widehat F(t))$ a set of all (classes of the equivalence which contain) polynomials of
type~(\ref{f2_12})
from $L_2(\mathcal{H}; [0,2\pi]; d\widehat F(t))$.
Set $L_{2;0}(\mathcal{H}; [0,2\pi]; d\widehat F(t)) = \overline{P (\mathcal{H}; [0,2\pi]; d\widehat F(t))}$.
Consider the following transformation:
$$ W_0 \left[
\sum_{j=0}^d e^{ijt} h_j 
\right]
= \sum_{j=0}^d x_{h_j,j},\qquad h_j\in \mathcal{H}, $$
which maps $P (\mathcal{H}; [0,2\pi]; d\widehat F(t))$ on the whole $\mathfrak{L}(\subseteq H)$.
Let us check that $W_0$ is well-defined. In fact, suppose that $p,q$ from~(\ref{f2_12})
belong to the same class of the equivalence. Then
$$ 0 = \left\|
\sum_{j=0}^d e^{ijt} (h_j-g_j) 
\right\|_{L_2}
=
\left(
\sum_{j=0}^d e^{ijt} (h_j-g_j),
\sum_{k=0}^d e^{ikt} (h_k-g_k) 
\right)_{L_2} = $$
$$ = \left( \sum_{j=0}^d x_{h_j-g_j,j}, \sum_{k=0}^d x_{h_k-g_k,k} \right)_H 
= \left\|
\sum_{j=0}^d x_{h_j,j} - \sum_{j=0}^d x_{g_j,j}
\right\|^2_H.
$$ 
Therefore $W_0$ is well-defined.
Moreover, $W_0$ is linear and relation~(\ref{f2_13}) shows that $W_0$ is isometric. By the continuity we extend $W_0$ to a unitary
transformation $W$ which maps $L_{2;0}(\mathcal{H}; [0,2\pi]; d\widehat F(t))$ on the whole $H$.
Denote by $U_0$ an operator in  $L_2(\mathcal{H}; [0,2\pi]; d\widehat F(t))$ which put into correspondence for a class of the equivalence
which has a representative 
$f(t)\in C_{00} (\mathcal{H};[0,2\pi])$  the class $[e^{it} f(t)]$.
It is readily checked that this operator is well-defined, linear and isometric. We extend it by the continuity
to a unitary operator $U$ on $L_2(\mathcal{H}; [0,2\pi]; d\widehat F(t))$.
Observe that
$$ A_0 x_{h,j} = x_{h,j+1} = W \left[ e^{i(j+1)t} h \right] = WU \left[ e^{ijt} h \right] = WUW^{-1} x_{h,j}, $$
for arbitrary $h\in \mathcal{H}$, $0\leq j\leq d-1$. Therefore
\begin{equation}
\label{f2_14}
WUW^{-1}\supseteq A.
\end{equation}
Set
$$ L_{2;1}(\mathcal{H}; [0,2\pi]; d\widehat F(t)) = 
L_2(\mathcal{H}; [0,2\pi]; d\widehat F(t))\ominus L_{2;0}(\mathcal{H}; [0,2\pi]; d\widehat F(t)), $$
and $\mathbf W = W\oplus E_{L_{2;1}(\mathcal{H}; [0,2\pi]; d\widehat F(t))}$.
Notice that $\mathbf W$ is a unitary transformation which maps $L_2(\mathcal{H}; [0,2\pi]; d\widehat F(t))$
on $H_1 := H\oplus L_{2;1}(\mathcal{H}; [0,2\pi]; d\widehat F(t))$.
Set $\widetilde A := \mathbf W U \mathbf{W}^{-1}$. By~(\ref{f2_14}) we see that $\widetilde A$ is a unitary
extension of $A$.
Let $\{ \widetilde E_t \}_{t\in [0,2\pi]}$ be a strongly left-continuous resolution of the identity of $\widetilde A$.
Denote by $\mathbf E_t$ ($t\in [0,2\pi]$) the corresponding strongly left-continuous spectral function of $A$.
For arbitrary $h,g\in \mathcal{H}$, $\zeta\in\mathbb{T}_e$ we may write 
$$ \int_0^{2\pi} \frac{1}{1-\zeta e^{it}} d(I^* \mathbf{E}_t I h,g)_H = 
\int_0^{2\pi} \frac{1}{1-\zeta e^{it}} d(\widetilde E_t x_{h,0}, x_{g,0})_{H_1} = $$
$$ =
\left( 
\left( E_{H_1} - \zeta\widetilde A 
\right)^{-1} \mathbf W [h], \mathbf W [g] 
\right)_{H_1} = 
\left( 
\mathbf W^{-1} \left( E_{H_1} - \zeta\widetilde A 
\right)^{-1} \mathbf W [h], [g] 
\right)_{L_2} =
$$
$$ = \left( 
\left( E_{L_2} - \zeta U 
\right)^{-1} [h], [g] 
\right)_{L_2} =
\int_0^{2\pi} \left(
d\widehat F(t) \frac{1}{1-\zeta e^{it}} h, g
\right)_{\mathcal{H}} =
$$
$$ = \int_0^{2\pi} 
\frac{1}{1-\zeta e^{it}} d(\widehat F(t)h, g)_{\mathcal{H}}. $$
Therefore 
$$ \int_0^{2\pi} \frac{1 + \zeta e^{it}}{1-\zeta e^{it}} d(I^* \mathbf{E}_t I h,g)_H = 
2 \int_0^{2\pi} \frac{1}{1-\zeta e^{it}} d(I^* \mathbf{E}_t I h,g)_H - $$
$$ - \int_0^{2\pi} d(I^* \mathbf{E}_t I h,g)_H = 
2 \int_0^{2\pi} \frac{1}{1-\zeta e^{it}} d(I^* \mathbf{E}_t I h,g)_H - (S_0 h, g)_{\mathcal{H}} = $$
$$ =
2 \int_0^{2\pi} 
\frac{1}{1-\zeta e^{it}} d(\widehat F(t)h, g)_{\mathcal{H}} - (\widehat F(2\pi) h,g)_{\mathcal{H}} = 
\int_0^{2\pi} 
\frac{1 + \zeta e^{it}}{1-\zeta e^{it}} d(\widehat F(t)h, g)_{\mathcal{H}}. $$
By the well-known inversion formula we conclude that $\widehat F(t) = I^* \mathbf E_t I$.

\begin{thm}
\label{t2_1}
Let the truncated operator trigonometric moment problem~(\ref{f1_1}) with $d\in \mathbb{N}$ be given and
condition~(\ref{f1_2}) hold. Let an operator $A_0$ in a Hilbert space $H$ be constructed as in~(\ref{f2_6}), 
$A = \overline{A_0}$.
All solutions of the moment problem have the following form:
\begin{equation}
\label{f2_15}
F(t) = I^* \mathbf{E}_t I,\qquad t\in [0,2\pi],
\end{equation}
where $I$ is defined by~(\ref{f2_8}) and $\mathbf{E}_t$ is a strongly left-continuous spectral function of $A$.
On the other hand, each strongly left-continuous spectral function of $A$ generates by~(\ref{f2_15})
a solution of the moment problem.
Moreover, the correspondence between all strongly left-continuous spectral functions of $A$ and
all solutions of the moment problem~(\ref{f1_1}) is one-to-one.
\end{thm}
\textbf{Proof.}
It remains to check that different left-continuous spectral functions of $A$ generate
different solutions of the moment problem~(\ref{f1_1}).
Set
$L_0 := \{ x_{h,0} \}_{h\in \mathcal{H}}$.
Choose an arbitrary element $x\in\mathfrak{L}$, $x=\sum_{j=0}^d x_{h_j,j}$, $h_j\in\mathcal{H}$.
For arbitrary $\zeta\in\mathbb{T}_e\backslash\{ 0\}$ there exists the following representation:
\begin{equation}
\label{f2_16}
x = v + y,\qquad v\in M_\zeta(A),\ y\in L_0.
\end{equation}
Here $v$ and $y$ may depend on the choice of $\zeta$.
In fact, for an arbitrary element $u\in D(A_0)$, $u=\sum_{j=0}^{d-1} x_{g_j,j}$, $g_j\in\mathcal{H}$, we may write
$$ (E_H - \zeta A) u = \sum_{j=0}^{d-1} x_{g_j,j} - \zeta \sum_{j=0}^{d-1} x_{g_j,j+1} = 
\sum_{j=0}^{d-1} x_{g_j,j} - \zeta \sum_{j=1}^{d} x_{g_{j-1},j} = $$
\begin{equation}
\label{f2_17}
= x_{g_0,0} + \sum_{j=1}^{d-1} x_{g_j - \zeta g_{j-1},j} + x_{-\zeta g_{d-1},d}. 
\end{equation}
Consider the following system of equations:
\begin{equation}
\label{f2_19}
\left\{ 
\begin{array}{cc}
g_j - \zeta g_{j-1} = h_j, & 1\leq j\leq d-1 \\
-\zeta g_{d-1} = h_d\end{array}
\right..
\end{equation}
We can find $g_{d-1}$, then $g_{d-2}$, ..., $g_0$.
Consider $u$ with this choice of $g_j$ and set $v := (E_H - \zeta A) u\in M_\zeta(A)$.
By~(\ref{f2_17}),(\ref{f2_19}) we see that
$$ x-v = x_{h_0,0} - x_{g_0,0} = x_{h_0-g_0,0} =: y\in L_0. $$
Then relation~(\ref{f2_16}) holds.

\noindent
Suppose to the contrary that two different strongly left-continuous spectral functions
$\mathbf E_{j,t}$, $j=1,2$, generate the same solution of the moment problem:
$I^* \mathbf E_{1,t} I = I^* \mathbf E_{2,t} I$. For arbitrary $f,g\in \mathcal{H}$ we have:
$$ (\mathbf E_{1,t} x_{f,0}, x_{g,0})_H = (\mathbf E_{2,t} x_{f,0}, x_{g,0})_H. $$
Multiplying by $\frac{1}{1-\zeta e^{it}}$ and integrating we get
\begin{equation}
\label{f2_22}
(\mathbf R_{1,\zeta} x_{f,0}, x_{g,0})_H = (\mathbf R_{2,\zeta} x_{f,0}, x_{g,0})_H,\qquad f,g\in\mathcal{H},\ \zeta\in\mathbb{T}_e,
\end{equation}
where $\mathbf R_{j,\zeta}$ is a generalized resolvent corresponding to $\mathbf E_{j,t}$, $j=1,2$.
Let $\mathbf R_{j,\zeta}$ is generated by a unitary extension $\widetilde A_j$ of $A$ in a Hilbert space
$\widetilde H_j\supseteq H$, $j=1,2$. Since for arbitrary $f\in D(A),\zeta\in\mathbb{T}_e$ and $j=1,2$ we have
$$ \left( E_{\widetilde H_j} - \zeta \widetilde A_j \right)^{-1} \left(E_H-\zeta A \right) f =
\left( E_{\widetilde H_j} - \zeta \widetilde A_j \right)^{-1} \left( E_{\widetilde H_j} - \zeta \widetilde A_j \right) f = f\in H, $$
then
\begin{equation}
\label{f2_24}
\mathbf R_{1,\zeta} u = \mathbf R_{2,\zeta} u,\qquad u\in M_\zeta(A),\ \zeta\in\mathbb{T}_e.
\end{equation}
Choose an arbitrary $\zeta\in \mathbb{T}_e\backslash\{ 0\}$. We may write:
$$ (\mathbf R_{j,\zeta} x_{f,0}, u) = (x_{f,0}, \mathbf R_{j,\zeta}^* u) = 
\left( x_{f,0}, u -  \mathbf R_{j, \frac{1}{\overline{\zeta}}} u \right), $$
where $f\in\mathcal{H}$, $u\in M_{ \frac{1}{\overline{\zeta}} }(A)$, $j=1,2$. 
By~(\ref{f2_24}) we get
\begin{equation}
\label{f2_26}
(\mathbf R_{1,\zeta} x_{f,0}, u)  = (\mathbf R_{2,\zeta} x_{f,0}, u),\qquad u\in M_{ \frac{1}{\overline{\zeta}} }(A),\ f\in \mathcal{H},\ 
\zeta\in\mathbb{T}_e\backslash\{ 0 \}.
\end{equation}
Choose an arbitrary element $w\in\mathfrak{L}$ and $\zeta\in \mathbb{T}_e\backslash\{ 0\}$. By~(\ref{f2_16}) we may write:
$w = v + y$, where $v\in M_{ \frac{1}{\overline{\zeta}} }(A)$, $y\in L_0$. 
By~(\ref{f2_22}),(\ref{f2_26}) we get
$(\mathbf R_{1,\zeta} x_{f,0}, w)  = (\mathbf R_{2,\zeta} x_{f,0}, w)$, $f\in\mathcal{H}$. Therefore
\begin{equation}
\label{f2_28}
\mathbf R_{1,\zeta} x = \mathbf R_{2,\zeta} x,\qquad x\in L_0,\ \zeta\in\mathbb{T}_e\backslash\{ 0 \}.
\end{equation}
For an arbitrary $w\in\mathfrak{L}$ using~(\ref{f2_16}) we may write:
$w = v' + y'$, where $v'\in M_\zeta(A)$, $y'\in L_0$, $\zeta\in\mathbb{T}_e\backslash\{ 0 \}$.
By~(\ref{f2_24}),(\ref{f2_28}) we get
$\mathbf R_{1,\zeta} w = \mathbf R_{2,\zeta} w$, $\zeta\in\mathbb{T}_e\backslash\{ 0 \}$.
Therefore $\mathbf E_{1,t}=\mathbf E_{2,t}$. This contradiction completes the proof.
$\Box$

Notice that relation~(\ref{f2_15}) is equivalent to the following relation:
\begin{equation}
\label{f2_30}
(F(t)g,h)_{\mathcal{H}} = (\mathbf{E}_t x_{g,0}, x_{h,0}),\qquad g,h\in\mathcal{H}.
\end{equation} 
From the latter relation it follows that
\begin{equation}
\label{f2_32}
\int_0^{2\pi} \frac{1+\zeta e^{it}}{1-\zeta e^{it}} d(F(t)g,h)_{\mathcal{H}} = 2(\mathbf{R}_\zeta x_{g,0}, x_{h,0})
-(S_0 g,h)_{\mathcal{H}},\quad g,h\in\mathcal{H},\ 
\zeta\in \mathbb{T}_e.
\end{equation}
By virtue of Chumakin's formula for the generalized resolvents of an isometric operator (see~\cite{cit_750_Ch}) we conclude that
the following formula:
$$ \int_0^{2\pi} \frac{1+\zeta e^{it}}{1-\zeta e^{it}} d(F(t)g,h)_{\mathcal{H}} = 
2 \left( 
\left[
E_H - \zeta (A\oplus\Phi_\zeta)
\right]^{-1} x_{g,0}, x_{h,0}
\right)_H -$$
\begin{equation}
\label{f2_34}
-
(S_0 g,h)_{\mathcal{H}},\qquad g,h\in\mathcal{H},\ 
\zeta\in \mathbb{D}.
\end{equation}
establishes a one-to-one correspondence between all functions $\Phi\in S(\mathbb{D};H\ominus D(A), H\ominus R(A))$ and
all solutions of the moment problem~(\ref{f1_1}).

We shall need the following proposition which is close to Frobenius's inversion formula for matrices.

\begin{prop}
\label{p2_1}
Let $M$ be a linear bounded operator in a Hilbert space $\mathbf{H}$, $D(M) = \mathbf{H}$. Suppose that
\begin{equation}
\label{fffffff2_53_5_1}
\mathbf{H} = \mathbf{H}_1 \oplus \mathbf{H}_2,
\end{equation}
where $\mathbf{H}_1$, $\mathbf{H}_2$ are subspaces of $\mathbf{H}$, and the operator $M$ has the following block representation:
\begin{equation}
\label{fffffff2_53_5_2}
M = \left(
\begin{array}{cc}
\mathbf{A} & \mathbf{B} \\
\mathbf{C} & \mathbf{D} \end{array}
\right),
\end{equation}
where $\mathbf{A} = P_{\mathbf{H}_1} M|_{\mathbf{H}_1}$,
$\mathbf{B} = P_{\mathbf{H}_1} M|_{\mathbf{H}_2}$,
$\mathbf{C} = P_{\mathbf{H}_2} M|_{\mathbf{H}_1}$, 
$\mathbf{D} = P_{\mathbf{H}_2} M|_{\mathbf{H}_2}$.
Suppose that $\mathbf{A}$ has a bounded inverse which is defined on the whole $\mathbf{H}_1$.
Then the following assertions hold.

\begin{itemize}
\item[(i)] If the operator $\mathcal{H} = \mathbf{D} - \mathbf{C} \mathbf{A}^{-1} \mathbf{B}$ has a bounded inverse which is defined on the whole 
$\mathbf{H}_2$, then
$M$ has a bounded inverse which is defined on the whole $\mathbf{H}$, and $M^{-1}$ has the following block representation
with respect to decomposition~(\ref{fffffff2_53_5_1}):
\begin{equation}
\label{fffffff2_53_5_3}
M^{-1} = \left(
\begin{array}{cc}
\mathbf{A}^{-1} + \mathbf{A}^{-1} \mathbf{B} \mathcal{H}^{-1} \mathbf{C} \mathbf{A}^{-1} & - \mathbf{A}^{-1} \mathbf{B} \mathcal{H}^{-1} \\
- \mathcal{H}^{-1} \mathbf{C} \mathbf{A}^{-1} & \mathcal{H}^{-1} \end{array}
\right).
\end{equation}

\item[(ii)] If the operator $M$ has a bounded inverse which is defined on the whole $\mathbf{H}$, then
the operator $\mathcal{H} = \mathbf{D} - \mathbf{C} \mathbf{A}^{-1} \mathbf{B}$ has a bounded inverse which is defined on the whole  $\mathbf{H}_2$,
and $M^{-1}$ has the block representation~(\ref{fffffff2_53_5_3}) with respect to decomposition~(\ref{fffffff2_53_5_1}).
\end{itemize}
\end{prop}

\textbf{Proof.}
$(i)$:
in this case the operator which is defined by the block representation in~(\ref{fffffff2_53_5_3}), is defined on the whole
$\mathbf{H}$ and it is bounded. The fact that it is the inverse of $M$ is verified by a direct block multiplication.

\noindent
$(ii)$: let us check that $\mathcal{H}$ is invertible. To the contrary, suppose that there exists a non-zero element
$h$ in $\mathbf{H}_2$: $\mathcal{H} h = 0$. Denote $u = - \mathbf{A}^{-1} \mathbf{B} h$. Then
\begin{equation}
\label{fffffff2_53_5_4}
\left(
\begin{array}{cc}
\mathbf{A} & \mathbf{B} \\
0 & \mathcal{H} \end{array}
\right)
\left(
\begin{array}{cc}
u \\
h \end{array}
\right) = 0.
\end{equation}
Notice that
$$ \left(
\begin{array}{cc}
\mathbf{A} & \mathbf{B} \\
0 & \mathcal{H} \end{array}
\right)
=
\left(
\begin{array}{cc}
E_{\mathbf{H}_1} & 0 \\
- \mathbf{C} \mathbf{A}^{-1} & E_{\mathbf{H}_2} \end{array}
\right)
\left(
\begin{array}{cc}
\mathbf{A} & \mathbf{B} \\
\mathbf{C} & \mathbf{D} \end{array}
\right),
$$
and the operator $\left(
\begin{array}{cc}
E_{\mathbf{H}_1} & 0 \\
- \mathbf{C} \mathbf{A}^{-1} & E_{\mathbf{H}_2} \end{array}
\right)$ has a bounded inverse which is defined on the whole  $\mathbf{H}$, equal to
$\left(
\begin{array}{cc}
E_{\mathbf{H}_1} & 0 \\
\mathbf{C} \mathbf{A}^{-1} & E_{\mathbf{H}_2} \end{array}
\right)$. 
Therefore an operator $\left(
\begin{array}{cc}
\mathbf{A} & \mathbf{B} \\
0 & \mathcal{H} \end{array}
\right)$ 
has a bounded inverse which is defined on the whole  $\mathbf{H}$. 
This contradicts to relation~(\ref{fffffff2_53_5_4}).
Thus, the operator $\mathcal{H}$ is invertible.
Since $R\left(
\begin{array}{cc}
\mathbf{A} & \mathbf{B} \\
0 & \mathcal{H} \end{array}
\right) = \mathbf{H}$, 
then for an arbitrary element $\widetilde h\in \mathbf{H}_2$ there exist elements $u_1\in \mathbf{H}_1$ и $h_2\in \mathbf{H}_2$, such that:
$$ \left(
\begin{array}{cc}
0  \\
\widetilde{h} \end{array}
\right) 
=
\left(
\begin{array}{cc}
\mathbf{A} & \mathbf{B} \\
0 & \mathcal{H} \end{array}
\right)
\left(
\begin{array}{cc}
u_1  \\
h_2 \end{array}
\right)
=
\left(
\begin{array}{cc}
\mathbf{A} u_1 + \mathbf{B} h_2  \\
\mathcal{H} h_2 \end{array}
\right). $$
Therefore $R(\mathcal{H}) = \mathbf{H}_2$. Since an operator $\mathcal{H}^{-1}$ is closed and defined on the whole $\mathbf{H}_2$, then
$\mathcal{H}^{-1}$ is bounded.
Applying assertion~$(i)$ we conclude that
$M^{-1}$ has a block representation~(\ref{fffffff2_53_5_3}) with respect to decomposition~(\ref{fffffff2_53_5_1}).
$\Box$

Return to our constructions for the solvable moment problem~(\ref{f1_1}) with $d\in \mathbb{N}$. 
Let us apply Proposition~\ref{p2_1} to the operator $M := E_H - \zeta (A\oplus\Phi_\zeta)$, $\zeta\in\mathbb{D}$,
$\Phi\in S(\mathbb{D};H\ominus D(A), H\ominus R(A))$,
in a Hilbert space $\mathbf{H} = H$ with $\mathbf{H_1} = D(A)$, $\mathbf{H_2} = H\ominus D(A)$.
In this case the operators $\mathbf{A}$, $\mathbf{B}$, $\mathbf{C}$, $\mathbf{D}$ in~(\ref{fffffff2_53_5_2})
have the following form:
$$ \mathbf{A} = E_{D(A)} - \zeta P^H_{D(A)} A,\quad
\mathbf{B} = -\zeta P^H_{D(A)} \Phi_\zeta,\quad $$
\begin{equation}
\label{f2_54}
\mathbf{C} = - \zeta P^H_{H\ominus D(A)} A,\quad
\mathbf{D} = E_{H\ominus D(A)} - \zeta P^H_{H\ominus D(A)} \Phi_\zeta,\qquad \zeta\in\mathbb{D}. 
\end{equation}
By~(\ref{f2_34}) we conclude that the following relation
$$ \int_0^{2\pi} \frac{1+\zeta e^{it}}{1-\zeta e^{it}} d(F(t)g,h)_{\mathcal{H}} = $$
\begin{equation}
\label{f2_56}
= 2
\left( 
\left\{
\mathbf{A}^{-1} + \mathbf{A}^{-1} \mathbf{B} \mathcal{H}^{-1} \mathbf{C} \mathbf{A}^{-1}
\right\} x_{g,0}, x_{h,0}
\right)_H
-
(S_0 g,h)_{\mathcal{H}},\qquad g,h\in\mathcal{H},\ 
\zeta\in \mathbb{D},
\end{equation}
establishes a one-to-one correspondence between all functions $\Phi\in S(\mathbb{D};H\ominus D(A), H\ominus R(A))$ and
all solutions of the moment problem~(\ref{f1_1}).
By the substitution of expressions from~(\ref{f2_54}) we obtain that relation~(\ref{f2_56}) takes the
following form:
$$ \int_0^{2\pi} \frac{1+\zeta e^{it}}{1-\zeta e^{it}} d(F(t)g,h)_{\mathcal{H}} = $$
$$ =
2
\left( 
\left\{
\mathbf{A}(\zeta) + \mathbf{B}(\zeta) \Phi_\zeta
\left(
E_{H\ominus D(A)} + \mathbf{C}(\zeta) \Phi_\zeta
\right)^{-1} 
\mathbf{D}(\zeta)
\right\} x_{g,0}, x_{h,0}
\right)_H
-
$$
\begin{equation}
\label{f2_58}
- (S_0 g,h)_{\mathcal{H}},\qquad
g,h\in\mathcal{H},\ \zeta\in \mathbb{D},
\end{equation}
where
$$ \mathbf{A}(\zeta) = 
\left(
E_{D(A)} - \zeta P^H_{D(A)} A
\right)^{-1},\quad
\mathbf{B}(\zeta) = -\zeta \mathbf{A}(\zeta) P^H_{D(A)}, $$
$$ \mathbf{D}(\zeta) = -\zeta P^H_{H\ominus D(A)} A \mathbf{A}(\zeta),\quad
$$
\begin{equation}
\label{f2_60}
\mathbf{C}(\zeta) = -\zeta P^H_{H\ominus D(A)} +
\zeta \mathbf{D}(\zeta) P^H_{D(A)},\qquad
\zeta\in \mathbb{D}.
\end{equation}
Finally, we obtain that
the following relation
$$ \int_0^{2\pi} \frac{1+\zeta e^{it}}{1-\zeta e^{it}} dF(t) = $$
\begin{equation}
\label{f2_62}
=
\mathcal{A}(\zeta) + \mathcal{B}(\zeta) \Phi_\zeta
\left(
E_{H\ominus D(A)} + \mathcal{C}(\zeta) \Phi_\zeta
\right)^{-1} 
\mathcal{D}(\zeta),\qquad 
\zeta\in \mathbb{D},
\end{equation}
establishes a one-to-one correspondence between all functions $\Phi\in S(\mathbb{D};H\ominus D(A), H\ominus R(A))$ and
all solutions of the moment problem~(\ref{f1_1}).
Here
$$ \mathcal{A}(\zeta) = 2 I^* \mathbf{A}(\zeta) I - S_0,\quad \mathcal{B}(\zeta) = 2 I^* \mathbf{B}(\zeta)|_{H\ominus R(A)}, $$
\begin{equation}
\label{f2_65}
\mathcal{C}(\zeta) = \mathbf{C}(\zeta)|_{H\ominus R(A)},\quad \mathcal{D}(\zeta) = \mathbf{D}(\zeta) I,\qquad
\zeta\in \mathbb{D}.
\end{equation}

\begin{thm}
\label{t2_2}
Let the truncated operator trigonometric moment problem~(\ref{f1_1}) with $d\in \mathbb{N}$ be given and
condition~(\ref{f1_2}) hold. Let an operator $A_0$ in a Hilbert space $H$ be constructed as in~(\ref{f2_6}), 
$A = \overline{A_0}$.
Relation~(\ref{f2_62}) establishes a one-to-one correspondence between all functions $\Phi\in S(\mathbb{D};H\ominus D(A), H\ominus R(A))$ and
all solutions of the moment problem~(\ref{f1_1}).
\end{thm}
\textbf{Proof.}
The proof follows from the preceding considerations. 
$\Box$

\begin{cor}
\label{c2_1}
Let the truncated operator trigonometric moment problem~(\ref{f1_1}) with $d\in \mathbb{N}$ be given and
condition~(\ref{f1_2}) hold. Let an operator $A_0$ in a Hilbert space $H$ be constructed as in~(\ref{f2_6}), 
$A = \overline{A_0}$.
The moment problem is determinate if and only if (at least) one of the defect numbers of $A$ is zero.
\end{cor}
\textbf{Proof.}
If one of the defect numbers of $A$ is zero, then the only function in $S(\mathbb{D};H\ominus D(A), H\ominus R(A))$ is zero.

On the other hand, if both defect numbers of $A$ are non-zero than we may choose unit nonzero vectors
$h\in H\ominus D(A)$, $g\in H\ominus R(A)$. Set $\Phi_1(\zeta)=0$, $\Phi_2(\zeta) = (\cdot, h) g$, $\zeta\in\mathbb{D}$.
Functions $\Phi_1,\Phi_2$ generate different solutions of the moment problem.
$\Box$

Let the truncated operator trigonometric moment problem~(\ref{f1_1}) with $d\in \mathbb{N}$ be given. 
Suppose that condition~(\ref{f1_2}) holds and \textit{the Hilbert space $\mathcal{H}$ is separable, $\mathcal{H}\not=\{ 0 \}$}.
Let $\mathfrak{C} = \{ f_k \}_{k=0}^{\omega-1}$, $1\leq \omega\leq +\infty$, is an orthonormal basis in $\mathcal{H}$.
Let us calculate the matrices of operators $\mathcal{A}$, $\mathcal{B}$, $\mathcal{C}$, $\mathcal{D}$ in~(\ref{f2_65})
with respect to some proper bases.
As a consequence, we can obtain the matrix of the operator appearing on the right of~(\ref{f2_62}) with respect to $\mathfrak{C}$
using the prescribed moments.

Observe that $H$ is separable. In fact, an arbitrary element $x$ of $\mathfrak{L}$ has the following form:
$x = \sum_{j=0}^d x_{h_j,j}$, $h_j\in \mathcal{H}$.
Let $\mathcal{W}$ be a dense subset of $\mathcal{H}$ (which is not supposed to be countable). 
Choose an arbitrary $\varepsilon > 0$.
There
exist elements $y_j\in \mathcal{W}$, $0\leq j\leq d$, such that
$$ \| h_j - y_j \|_{\mathcal{H}} \leq \frac{\varepsilon}{(d+1)\| S_0 \|^{ \frac{1}{2}} },\qquad 0\leq j\leq d. $$
Then
$$ \left\| x - \sum_{j=0}^d x_{y_j,j} \right\|_H \leq \sum_{j=0}^d \| x_{h_j,j} - x_{y_j,j} \|_H = $$
$$ = \sum_{j=0}^d \sqrt{ ( x_{h_j-y_j,j}, x_{h_j-y_j,j} ) } = \sum_{j=0}^d 
\sqrt{ ( S_0 (h_j-y_j), h_j-y_j )_{\mathcal{H}} } \leq  $$ 
$$ \leq \sum_{j=0}^d \| S_0 \|^{\frac{1}{2}}
\| h_j-y_j \|_{\mathcal{H}} \leq \varepsilon. $$ 
Therefore a set $\widetilde{\mathcal{W}} := 
\left\{
\sum_{j=0}^d x_{y_j,j},\quad y_j\in \mathcal{W}
\right\}
$ is dense in $H$. If $\mathcal{W}$ is chosen countable then $\widetilde{\mathcal{W}}$ is countable, as well.
Thus, $H$ is separable. On the other hand, if we choose $\mathcal{W} = \mathop{\rm Lin}\nolimits \mathfrak{C}$ then
we obtain that $\mathop{\rm Lin}\nolimits\{ x_{f_k,j} \}_{k\in\overline{0,\omega - 1},\ 0\leq j\leq d}$ is dense in $H$.

At first we suppose that \textit{the moment problem is indeterminate}.
If we would have $D(A) = \{ 0 \}$, then $(S_0 h,g)_{\mathcal{H}} = (x_{h,0}, x_{g,0})_H = 0$,
$h,g\in\mathcal{H}$. This would meant that $S_0 = 0$ and the moment problem is determinate.
This contradicts to our assumptions. Therefore $D(A)\not= \{ 0 \}$.
Observe that 
\begin{equation}
\label{f2_66}
D(A) = \mathop{\rm span}\nolimits\{ x_{f_k,j} \}_{k\in\overline{0,\omega - 1},\ 0\leq j\leq d-1}. 
\end{equation}
Let us numerate the elements of the set $\Omega := \{ x_{f_k,j} \}_{k\in\overline{0,\omega - 1},\ 0\leq j\leq d-1}$
by a unique index: $\Omega := \{ y_n \}_{n\in\overline{0,\kappa-1}}$, $0\leq\kappa\leq +\infty$. 
Apply the Gram-Schmidt orthogonalization procedure to the following sequence:
\begin{equation}
\label{f2_67}
\{ y_n \}_{n\in\overline{0,\kappa-1}},
\end{equation}
removing linear dependent elements, if they appear.
After the orthogonalization we shall
obtain an orthonormal basis $\mathfrak{A}_1 = \{ u_{k} \}_{k=0}^{\tau-1}$, $1\leq\tau\leq \infty$
in $D(A)$.
Apply the Gram-Schmidt orthogonalization procedure to the following sequence:
\begin{equation}
\label{f2_68}
\left\{ x_{f_k,d} - P^H_{D(A)} x_{f_k,d} \right\}_{k\in\overline{0,\omega - 1}}, 
\end{equation}
removing linear dependent elements, if they appear.
Observe that the elements $P^H_{D(A)} x_{f_k,d}$ in~(\ref{f2_68}) can be constructed using $\mathfrak{A}_1$.
After this orthogonalization we shall
obtain an orthonormal basis $\mathfrak{A}_2 :=  \{ u_j' \}_{j=0}^{\delta-1}$, $1\leq\delta\leq +\infty$ in $H\ominus D(A)$.

Observe that
$\mathfrak{A}_1' := \{ A_0 u_k \}_{k=0}^{\tau-1}$ is an orthonormal basis in $R(A)$.
Apply the Gram-Schmidt orthogonalization procedure to the following sequence:
\begin{equation}
\label{f2_69}
\left\{ x_{f_k,0} - P^H_{R(A)} x_{f_k,0} \right\}_{k\in\overline{0,\omega - 1}},
\end{equation}
removing linear dependent elements, if they appear.
Observe that the elements $P^H_{R(A)} x_{f_k,0}$ in~(\ref{f2_68}) can be constructed using $\mathfrak{A}_1'$.
By the orthogonalization 
we shall obtain 
$\mathfrak{A}_2' := \{ v_k' \}_{k=0}^{\rho -1}$, $1\leq \rho\leq +\infty$.
Notice that $\mathfrak{A}_2'$ is an orthonormal basis in $H\ominus R(A)$.

The above orthonormal bases can be used to construct the matrices of operators on the right in~(\ref{f2_62}).
Here relation~(\ref{f2_5}) will be used intensively.

In the case of \textit{the determinate moment problem} the right-hand side of~(\ref{f2_62}) is
equal to $\mathcal{A}(\zeta)$. Thus, in this case we can use $\mathfrak{C}$ and an orthonormal basis $\mathfrak{A}$,
constructed by the orthogonalization of $\{ x_{f_k,j} \}_{k\in\overline{0,\omega - 1},\ 0\leq j\leq d}$.

\begin{center}
{\large\bf On the truncated operator trigonometric moment problem.}
\end{center}
\begin{center}
{\bf S.M. Zagorodnyuk}
\end{center}

In this paper we study the truncated operator trigonometric moment problem. All solutions of the moment problem
are described by a Nevanlinna-type parameterization. In the case of moments acting in a separable Hilbert space, the matrices of
the operator coefficients in the Nevanlinna-type formula are calculated by the prescribed moments.
Conditions for the determinacy of the moment problem are given, as well.
}


\end{document}